\documentclass[smallextended]{svjour3}       
\smartqed  
\usepackage{graphicx}
\usepackage{booktabs}
\usepackage{amsmath}
\usepackage{amsfonts}
\usepackage{amssymb}
\usepackage{epstopdf} 
\usepackage{bm}
\usepackage{color}
\usepackage{mathptmx}
\usepackage{mathrsfs}  
 
 \newtheorem{assumption}{Assumption}
\newcommand{\x}{\mathbf{x}}
\newcommand{\s}{\mathbf{s}}

\newcommand{\vertiii}[1]{{\left\vert\kern-0.25ex\left\vert\kern-0.25ex\left\vert #1 
		\right\vert\kern-0.25ex\right\vert\kern-0.25ex\right\vert}}
		
\def\eop{{\ \vrule height 5pt width 5pt depth 0pt}}

\begin{document}

\title{An Operator-Splitting  Finite Element Method for the Numerical Solution  of Radiative Transfer Equation
}

\titlerunning{An  Operator-Splitting FEM for Numerical Solution of  RTE}        

\author{Sashikumaar Ganesan        \and
        Maneesh Kumar Singh 
}


\institute{S. Ganesan   \at
              Department of Computational and Data Sciences,\\ Indian Institute of Science, Bangalore, India, 560012 \\
              \email{sashi@iisc.ac.in}           
           \and
           M. K. Singh \at
           Department of Computational and Data Sciences,\\ Indian Institute of Science, Bangalore, India, 560012\\
           \email{maneeshsingh@iisc.ac.in, maneesh.sr51@gmail.com} 
}

\date{Received: date / Accepted: date}

\maketitle

\begin{abstract}
An operator-splitting finite element scheme for the time-dependent, high-dimensional radiative transfer equation is presented in this paper. 
The streamline upwind Petrov-Galerkin finite element method and discontinuous Galerkin finite element method are used for the spatial-angular discretization of the radiative transfer equation, whereas the implicit backward Euler scheme is used for temporal discretization. 
Error analysis of the proposed numerical scheme for the fully discrete radiative transfer equation is presented. 
The stability and convergence estimates for the fully discrete problem are derived.
Moreover, an operator-splitting algorithm for numerical simulation of high-dimensional equations is also presented. The validity of the derived estimates and implementation is illustrated with suitable numerical experiments.

\keywords{Radiative transfer equation \and Operator-splitting method \and  Streamline upwind Petrov Galerkin finite element methods \and Backward Euler scheme \and Stability and Convergence analysis  }
\subclass{MSC 65M12 \and 65M15 \and 65M60  \and 65R20}
\end{abstract}

\section{Introduction}
Radiation plays a significant role, both as a detectable and dominant mechanism for transmitting energy inside and outside a system, in several areas, including optics, astrophysics,  atmospheric science, and remote sensing.
Consequently, the propagation of radiation through a medium is one of the most critical processes studied extensively. Analyzing the released radiation from an object provides insight into the radiative source, the medium between the object and the observer, and its surroundings. Modeling this physical process results in a time-dependent, six-dimensional Partial Differential Equation (PDE). The higher dimension of PDE is one of the challenges associated with the solution of the radiative transfer equation (RTE).
Besides, considerable uncertainty is added due to radiation's ability to
affect the medium's state, which is the source of the radiation itself. Although analytic solutions to RTE exist for simple cases, numerical solutions are often sought for more realistic, complex applications. 
Therefore, it is exciting and, at the same time challenging to develop numerical schemes for the radiative model. More details on the radiative transfer model can be found in \cite{chandrasekhar,kanschatnumerical}. The design and implementation of the computational method for time-dependent, high-dimensional radiative transfer equations remain a challenging task in computational science, even though tremendous advances have been made in this area over the past few years. 

\subsection{Model problem}
Let $\Omega$ be a bounded domain in $\mathbb{R}^{3}$  with a smooth boundary $\partial \Omega$. Denote by $\mathbf{n}(\mathbf{x})$
the unit outward normal for $\x \in \partial \Omega$. Let the angular space $S^{2}$ be  the unit sphere in $\mathbb{R}^{3}$. For each fixed direction $\s(s_1, s_2, s_3) \in S^{2}$, we introduce the following
subsets of the boundary $\partial \Omega$:
\[
\partial \Omega_{\s,-}=\{\x \in \partial \Omega: \s \cdot \mathbf{n}(\x)<0\},\quad \partial \Omega_{\s,+}=\{\x \in \partial \Omega: \s \cdot \mathbf{n}(\x) \geq 0\}.
\]
Let $\mathcal{G}$ be the whole domain which is the tensor product
domain in space and angle. Then the boundary $\Gamma=\partial \Omega \times S^{2}$ can be split into two parts
\[
\Gamma_{-}=\{(\x,\s): \x \in \partial \Omega_{\s,-},\,\s \in S^{2}\},\quad 
\Gamma_{+}=\{(\x,\s): \x \in \partial \Omega_{\s,+},\,\s \in S^{2}\}
\]
as the inlet and outlet boundaries.
In this article, the high-dimensional radiative transfer equation (RTE) is formally defined by 
an initial-boundary-value problem:
\begin{align} 
\displaystyle \dfrac{\partial u}{\partial t}+\s \cdot \nabla u+\sigma_{\tau} u-\sigma_{s} \int_{S^{2}}u (t,\x,\s')\Phi(\s, \s')d \s'  &=f,\, \mathrm{in}~ (0,T]\times\Omega \times S^{2}, \nonumber \\ 
u(t=0,\x,\s)&=u_{0},\, \mathrm{in}~ \Omega\times S^{2}, \label{RT1}\\ 
u(t,\x,\s)&=0,~\,\mathrm{on}~ (0,T]\times\Gamma_{-}, \nonumber
\end{align}
where 
$\sigma_{\tau}(\x)=\sigma_{a}(\x)+\sigma_{s}(\x)$. Here, $\sigma_{a}(\x)$ and $\sigma_s(\x)$ are the total absorption and scattering coefficients, respectively. 
For simplicity, the particle speed is assumed to be one.
Here, the scattering phase function  $\Phi(\s, \s')$  describes the probability of a photon at position $\x$ that originally propagates in the direction $\s$, and $\s'$ as its new propagation direction after the scattering event. Note that the angular variable $\s$ in the spherical coordinate system is denoted as 
$
\s=(\sin{\theta}\cos{\phi},\sin{\theta}\sin{\phi},\cos{\theta})^{T}.$
Also,  we make the  following assumptions on  the data of the model problem (\ref{RT1}) as
Further, the data $\sigma_{\tau},\sigma_{s},f$ and $u_{0}$ of the model problem (\ref{RT1}) are assumed to be sufficiently smooth. Note that,  the given RTE  model can  also be viewed as   high-dimensional integro-differential equation.

The existing numerical schemes for RTE models can be classified into (i)  stochastic approach and (ii) deterministic approach. Among all stochastic approaches, the Monte Carlo method is often used to solve the radiative transfer equation, see, for example,  \cite{gentile2001implicit,halton,howell} and the references therein.
Nevertheless, the Monte-Carlo simulation's computational cost is very high due to its iterative design, and it increases when the optical depth becomes large.

Several deterministic numerical schemes have been proposed in the literature for the stationary RTE, see for example, \cite{avila2011spatial,badri1,castro_JQSRT19,Egger_SIAM16,egger_SIAM1,ren2019fast} and the references therein.  A robust numerical $S_{N}$-DG-approximations for radiation transport has been discussed in \cite{kanschat_SIAMJSC14,ragusa2012robust}. Stabilized finite element scheme with discrete ordinate method has been discussed for steady-state RTE models in \cite{han,wang}. A numerical scheme based on Ad hoc angular discretization and vectorial finite elements for spatial discretization has been studied in \cite{Favennec_JQSRT19}. Recently, an adaptive nested source term iteration method for steady-state RTE has been presented in \cite{dahmen2020adaptive}. For a time-dependent RTE model in one-dimensional slab geometry, a semi-analytical numerical method has been presented in \cite{de2002semi}. A low-rank approximation for time-dependent radiation transport in one- and two-dimensional Cartesian geometries has been discussed in \cite{peng2020low}. In \cite{dehghanian2020transient},  a variable discrete ordinates method has been used to solve the transient radiation heat transfer in a semi-transparent slab. Despite several numerical schemes proposed in the literature, numerical solution of the time-dependent high-dimensional RTE is still challenging and is an active research field.

The operator-splitting finite element methods have been developed in the recent past for many high-dimensional physical and mathematical models. For example, an operator-splitting numerical method for the micro-macro dilute polymeric fluid model has been provided in \cite{knezevic_ESIAM}. For a high-dimensional
the convection-diffusion problem, an operator-splitting method with detailed numerical implementation, has been presented in \cite{ganeshan_AMC1}. The high-dimensional population balance equation using the operator-splitting method has been discussed in \cite{ankersashi_CCE,ganeshan_ESAIM1}. An operator-splitting finite element method for an efficient parallel solution of high-dimensional population balance systems has been discussed in \cite{ganeshan_CEC}. More details on the standard operator-splitting FEM can be found in book  \cite{ganesan2017finite}.

For high-dimensional time-dependent RTE, we present an operator-splitting heterogeneous FEM. Further, \emph{a priori} error estimate for the proposed numerical scheme is presented, which is the main contribution of this research article.
The key idea is to split the RTE model problem concerning the internal (angular) and the external (spatial) directions, resulting in a transient transport problem and a time-dependent integro-differential equation. The transient transport problem is numerically approximated using the streamline upwind Petrov-Galerkin (SUPG) finite element method, whereas the discontinuous Galerkin method with piecewise constant polynomials (DG($0$)) is used for the integro-differential equation.   The proposed numerical scheme uses this tailor-made spatial discretization method and the implicit backward Euler scheme. The stability estimate of the fully discrete form of the proposed scheme is first derived. A convergence analysis is then established under the assumption of certain regularity conditions on the data and the stabilization parameter. Finally, an array of numerical experiments is provided to support the theoretical error estimates of the numerical approximation.

The rest of the article is organized as follows. In Section \ref{sec2}, we briefly discuss the weak formulation of the operator-splitting method for the model problem, and later, we discuss the finite element approximation of RTE. Further, a fully discrete form is derived in this section. Next,  the stability estimate of the discrete problem and the convergence analysis of the numerical approximation is presented in Section \ref{sec3}. Further,  the implementation of the numerical scheme is discussed in Section \ref{sec4}. Finally, a concluding remark is discussed in Section \ref{sec5}.

\section{Finite element approximation for RTE}\label{sec2}
This section starts with   prerequisites for the finite element discretization of  the model problem (\ref{RT1}). 
Let $L^{2}$ and $H^{m}$ be the Sobolev spaces.  Denote the $L^{2}$-inner product with respect to spatial variable $\x$ over the domain $\Omega$ as $(\cdot,\cdot)_{\x}$.  
The inner product  and  $L^{2}$-norm over the entire domain  $\Omega \times S^{2}$  are defined  by
\[
\begin{array}{ll}
\displaystyle (v,w) := \int_{S^{2}}(v,w)_{\x}\,d \s=\int_{S^{2}}\int_{\Omega}v w \, d\x \,d \s,\qquad \|v\|^{2}_{0} =(v,v), \\[8pt]
\displaystyle  \langle v, w \rangle_{\Gamma_+}  := \int_{S^{2}}\int_{\partial \Omega_{\s,+}}(v,w)_{\x}\, d\x \, d \s, \qquad \|v\|^{2}_{\Gamma_+} =\langle v, v \rangle_{\Gamma_+}.
\end{array}
\]
For simplification of mathematical presentation, we have  omitted $ d\s$ from $\int_{S^{2}}(v,w)_{\x}\,d \s$ and we simply write $\int_{S^{2}}(v,w)_{\x}$ throughout this article. 
We will adopt the notation associated with the operator-splitting technique introduced in \cite{ganeshan_AMC1,ahamd_ANM} for the numerical analysis of this article. 

Next, we introduce  Bochner spaces.  
Let $Z$ be a Banach space associated with the spatial variable $\x$ equipped with the norm $\|\cdot\|_Z$. For spaces $Z$ and $Y$,
we use a short notation $Y(Z) := Y (S^{2}; Z)$ and   define the following spaces
\[
\begin{array}{ll}
\mathbb{C}(S^{2}; Z):=\left\{v:\Omega\rightarrow Z \big| v \,\,\mbox{continuous},\,\,\displaystyle \sup_{\s \in S^{2}}\|v(\s)\|_{Z} < \infty \right\},\\[12pt]
L^{2}(S^{2}; Z):=\left\{v:\Omega\rightarrow Z \bigg| \displaystyle \int_{S^{2}}\|v\|^{2}_{Z} < \infty\right\},\\[12pt]
H^{m}(S^{2}; Z):=\left\{v \in L^{2}(S^{2}; Z)\bigg|\dfrac{\partial^{j}v}{\partial \s^{j}} \in L^{2}(S^{2}; Z),\,\, 1 \leq j \leq m\right\},
\end{array}
\]
where the  derivatives $\partial^{j}v/\partial \s^{j}\,$($j$ a multi-index)
are expressed in the sense of distributional derivative on $S^{2}$ and $m$ is an integer. The  norms in the above  spaces are given by
\[
\|v\|_{C(Z)}:=\sup_{\s \in S^{2}}\|v\|_{Z}, \qquad \|v\|^{2}_{L^{2}(Z)}:=\int_{S^{2}}\|v\|^{2}_{Z}, \quad 
\|v\|^{2}_{H^{m}(Z)}:=\int_{S^{2}}\sum_{|j|\leq m}\left\|\dfrac{\partial^{j}v}{\partial \s^{j}}\right\|^{2}_{Z}.
\]

Before deriving  the numerical approximation of the model problem (\ref{RT1}), the  analytical properties of the scattering phase function  is discussed here.
An operator $\mathcal{K}$ is defined  by 
\[
\mathcal{K}v(t,\x,\s)=\int_{S^{2}}\Phi(\s, \s')v(t,\x,\s')d \s'.
\]
\begin{assumption}\label{Asm1}
Assume that the scattering kernel $\Phi$ holds the following conditions:
\begin{itemize}
    \item $\Phi$ is a measurable function and positive, \emph{i.e.}, $\Phi(\s, \s') \geq 0$ for  $\s, \s' \in S^{2}$, 
    \item  and it satisfies 
    \begin{equation}\label{normelize}
\int_{S^{2}}\Phi(\s, \s')d \s'=1.	
\end{equation}
\end{itemize}
\end{assumption}

\begin{lemma}\label{prop_k1}
$\mathcal{K} : \Omega \times S^{2} \rightarrow \Omega \times S^{2}$ is a self-adjoint and bounded linear operator with
$\|\mathcal{K}v\|_{0} \leq \|v\|_{0}$, for all $v \in \Omega \times S^{2}$.
\end{lemma}
\textbf{ Proof.} The proof is given  in  \cite[Lemma 2.6]{egger2012mixed}.
\hfill \eop

\vskip 2mm
Next,  denote the  {\it removal  operator} \cite{egger2012mixed} by
\[
\mathcal{K}_{\sigma}v(t,\x,\s)=\sigma_{a}(\x)v(t,\x,\s)+\sigma_{\s}(\x)v(t,\x,\s)-\mathcal{K}v(t,\x,\s),
\]
where $\sigma_{a}(\x)v(t,\x,\s)$  models the absorption of particles
by the medium. The absorption and remission of particles during the scattering process is described by $(\sigma_{\s}(\x)v(t,\x,\s)-\mathcal{K}v(t,\x,\s))$.

\begin{assumption}\label{Asm2} The absorption and scattering coefficient satisfy following conditions 
\begin{itemize}
    \item $\sigma_{\s}$ is measurable, non-negative, and uniformly bounded, \emph{i.e.}, there exists   $\overline{\sigma}_{\s} \in \mathbb{R}^{+}$ such that $0 \leq \sigma_{\s}(\x) \leq  \overline{\sigma}_{\s} $ for \emph{a.e} $\x \in \Omega$.
    \item $\sigma_{a}$ is measurable, non-negative, and uniformly bounded, \emph{i.e.}, $0 < \underline{\sigma}_{a} \leq \sigma_{a}(\x) \leq  \overline{\sigma}_{a} $ for \emph{a.e} $\x \in \Omega$ and  $\underline{\sigma}_{a}, \overline{\sigma}_{a} \in \mathbb{R}^{+}$. For the convergence analysis, we assume that $\underline{\sigma}_{a} \geq 1/8$.
\end{itemize}
\end{assumption}

\begin{lemma}\label{prop_k2} The operator
$\mathcal{K}_{\sigma} : \Omega \times S^{2} \rightarrow \Omega \times S^{2}$ is a self-adjoint and elliptic bounded linear operator and it satisfies  following estimates:
\begin{equation}\label{scat_prop}
\begin{array}{ll}
    (\mathcal{K}_{\sigma} v,v)  & \geq \underline{\sigma}_{a} \|v\|^{2}_{0},\\[6pt] 
    (\mathcal{K}_{\sigma} v,w)  &  \leq \|\mathcal{K}_{\sigma} v\|_{0}\|w\|_{0} \leq (2\overline{\sigma}_{\s}+\overline{\sigma}_{a})  \| v\|_{0}\|w\|_{0}.
\end{array}
\end{equation}
\end{lemma}
\textbf{ Proof.} The proof  of lemma is discussed in  \cite[Lemma 2.7]{egger2012mixed}.
\hfill \eop

\subsection{Operator-Splitting Method}

The gradient operator in the RTE model (\ref{RT1}) is defined for the spatial variable $\x$ only. Thus, we can take advantage of the decomposition of the model problem by decomposing the model problem  (\ref{RT1}) into a purely convective problem in the space and an integro-differential equation in the angular variable.  Let $0 = t_{0}<t_{1}<\ldots<t_{N}=T$ be the time discretization of the time interval $[0,T]$.
Using the Lie's operator-splitting method, in the time
interval $(t_{n} , t_{n+1} )$, the operator-splitting method  of  the model problem (\ref{RT1}) read:

\vskip 2mm
\noindent
\textbf{Step 1. ($\s$-direction)}

For given $\tilde{u}(t^n)=u(t^n)$, find $\tilde{u}:(t^{n}, t^{n+1})\times \Omega \times S^{2}\rightarrow \mathbb{R}$  such that,
\begin{equation}\label{seq}
\begin{array}{ll}
\displaystyle \dfrac{\partial \tilde{u}}{\partial t}+\mathcal{K}_{\sigma}\tilde{u} =0, \,\,(t^{n}, t^{n+1})\times \Omega \times S^{2},\\[12pt]
\tilde{u}(t^{n},\x,\s)=u(t^{n},\x,\s),
\end{array}
\end{equation}
by considering $\x$ as a parameter. In this step, the solution is updated in the $\s$-direction. Then, this solution $ \tilde{u} $
is taken as the initial solution for the $\x$-direction update.

\vskip 2mm
\noindent
\textbf{Step 2. ($\x$-direction)}

For given $u(t^n)=\tilde{u}(t^{n+1})$, find $u:(t^{n}, t^{n+1})\times \Omega \times S^{2}\rightarrow \mathbb{R}$  such that 
\begin{equation}\label{xeq}
\begin{array}{ll}
\displaystyle \dfrac{\partial u}{\partial t}+\s \cdot \nabla u=f, \,\,(t^{n}, t^{n+1})\times \Omega \times S^{2},\\[8pt]
u=0,\,\,(t^{n}, t^{n+1})\times \Gamma_{-} ,\\[8pt]
u(t^{n},\x,\s)=\tilde{u}(t^{n+1},\x,\s).
\end{array}
\end{equation}
by considering the variable  $\s$ as a parameter. Here, the solution $u$ in the time step $(t^{n}, t^{n+1})$  is obtained by first updating in
$\s$-direction~\eqref{seq} and then updating in $\x$-direction~\eqref{xeq}.


A weak formulation for both   steps \eqref{seq} and \eqref{xeq}   will be introduced below. 
Let  $\tilde{V}$ is defined by
$$\tilde{V}=\{v| v \in L^{2}(\Omega),\,\s \cdot \nabla v \in L^{2}(\Omega),\,\,|\s \cdot \mathbf{n}|^{1/2}v \in L^{2}(\partial \Omega_{\s,\pm})\}.$$

We denote by $V=\{v \in \tilde{V}: v|_{\partial \Omega_{\s,-}}=0\}$ and   $W:=L^{2}(S^{2})$. Finally, we introduce 
\[
\mathcal{P}:=\left\{v \in  L^{2}(\Omega \times
S^{2})|v \in  L^{2}(\Omega ;
W) \cap v \in  L^{2}(S^{2} ;V)\right\}.
\]
From definition of the finite element space $\mathcal{P}$, any smooth function $v \in L^{2}(S^{2} ;V) \subset \mathcal{P}$ satisfies 
\[
 \displaystyle \int_{S^{2}} \|v\|^{2}_{L^{2}(\Omega)} < \infty, \quad \int_{S^{2}} \|\s \cdot \nabla v\|^{2}_{L^{2}(\Omega)}< \infty, \quad \int_{S^{2}} \||\s \cdot \mathbf{n}|^{1/2}v\|^{2}_{L^{2}(\partial \Omega_{\s,\pm})}< \infty.
\]

Now, we introduce the weak form  for the operator-splitting method (\ref{seq}) and (\ref{xeq}), which  is given by
\vskip 3mm
\noindent
\textbf{Step 1.} Find $\tilde{u}:(t^{n},t^{n+1})\rightarrow \mathcal{P}$ with $\tilde{u}(t^{n})=u(t^{n})$ such that
\begin{equation}\label{fe1}
\int_{S^{2}}(\tilde{u}_{t},v)_{\x}\,d\s+\int_{S^{2}}(\mathcal{K}_{\sigma}\tilde{u},v)_{\x}\,d\s=0,\,\,\forall v \in \mathcal{P},
\end{equation}
The norm associated with the weak formulation (\ref{fe1}) is simply the inner product norm $\|\cdot\|_{0}$.

\vskip 3mm
\noindent
\textbf{Step 2.} Find $u:(t^{n},t^{n+1})\rightarrow \mathcal{P}$ with $u(t^{n})=\tilde{u}(t^{n+1})$ such that
\begin{equation}\label{fe2}
\int_{S^{2}}(u_{t},v)_{\x}+\int_{S^{2}}a(u,v)=\int_{S^{2}}(f,v)_{\x},\,\,\forall v \in \mathcal{P},
\end{equation}
where the bilinear form  $a(u,v)=(\s \cdot \nabla u,v)_{\x}$. It is well documented in \cite{Eril_CMAME1}  that standard FEM are known to produce spurious oscillation.  To achieve the coercivity  of the  bilinear form $a(u,v)$, the test function is taken as $(v+\delta \s \cdot \nabla v)$, where $\delta$ is the stabilization parameter. A more detail about the $\delta$ is given later.  Finally,  the norm is given by
\[
\vertiii{v}^{2}= \int_{S^{2}} \left(\delta\|\s \cdot \nabla v\|^{2}_{L^{2}(\Omega)}+ \||\s \cdot \mathbf{n}|^{1/2}v\|^{2}_{L^{2}(\partial \Omega_{\s,\pm})}\right). 
\]


Here, we will briefly discuss the existence and uniqueness of the weak formulation (\ref{fe1}) and (\ref{fe2}) . We start with the weak form (\ref{fe1}). By following 
\cite[Theorem 6.1]{ganesan2017finite}, it is enough to show that 
\begin{align}
    \int_{S^{2}}(\mathcal{K}_{\sigma}\tilde{u},v)_{\x}\,d\s &\leq M_{1} \|\tilde{u}\|_{0} \|v\|_{0}, \,\, \mbox{a.e.}\, 0<t< T,\,\,\tilde{u}, v \in \mathcal{P},\label{eu1}\\[12pt]
     \int_{S^{2}}(\mathcal{K}_{\sigma}\tilde{u},\tilde{u})_{\x}\,d\s &\geq \alpha_{1}  \|\tilde{u}\|^{2}_{0}, \mbox{a.e.}\, 0<t< T,\,\,\tilde{u}\in \mathcal{P},\label{eu2}
\end{align}
where $M_{1}, \alpha_{1}$ are the positive constants. 

By using  Lemma \ref{prop_k2}, the required results (\ref{eu1}) and (\ref{eu2}) can be easily verified. Next, we  prove the existence and uniqueness of  (\ref{fe2}).  It is sufficient  to show that 
\begin{align}
    \int_{S^{2}}a(u,v+\delta \s \cdot \nabla v)\,d\s &\leq M_{2} \vertiii{u} \vertiii{v}, \,\, \mbox{a.e.}\, 0<t< T,\,\,\tilde{u}, v \in \mathcal{P},\label{eu3}\\[12pt]
     \int_{S^{2}}a(u,u+\delta \s \cdot \nabla u)\,d\s &\geq \alpha_{2}  \vertiii{u}^{2}, \mbox{a.e.}\, 0<t< T,\,\, u \in \mathcal{P},\label{eu4}
\end{align}
The inequality (\ref{eu3}) can be easily obtain by using Cauchy–Schwarz (C-S) inequality. And 
\[
\begin{array}{ll}
\displaystyle \int_{S^{2}}a(u,u+\delta \s \cdot \nabla u)
&= \displaystyle \int_{S^{2}} \int_{\partial \Omega_{s,+}}\frac{1}{2}(\s \cdot \mathbf{n}) u^{2}+\int_{S^{2}} \int_{\Omega} \delta (\s \cdot \nabla u )^{2}\\[16pt]
& \geq \dfrac{1}{2}\vertiii{u}^{2}_{2}, 
\end{array}
\]
which proves (\ref{eu4}) for $\alpha_{2} =1/2$. This completes the discussion of the existence and uniqueness of the weak formulation (\ref{fe1}) and (\ref{fe2}).

\subsection{Angular and spatial discretization}
In this current  subsection, we derive a semi-discrete
form of the operator-split equations.
It is well-known that the standard Galerkin finite element method for convection problems  (\ref{fe2}) induces spurious oscillation in the numerical solution. Therefore, we prefer the SUPG method for  spatial discretization. Since the numerical approximation of the double integral term in \eqref{fe1} will be compute-intensive, we implement DG($0$) for the angular discretization. 

Let $ S^{2}_{h}$ be a subdivision of $S^{2}$ into a surface mesh, which is obtained by discretizing the unit sphere $S^{2}$ using hierarchical sectioning of the sphere into spherical triangles. In particular, the   subdivision $ S^{2}_{h}$ is obtained by projecting polyhedra onto the unit sphere. 
More details on this type of triangulation can be found in  \cite[Chapter 3]{kanschat1996parallel} and \cite{kophazi2015space}. Further, the mesh size of the spherical triangles $K_{\s}$ in  $ S^{2}_{h}$ are denoted   by
\[ 
h_{\s}:= \max_{K_{\s} \in  S^{2}_{h}} h_{K_{\s}}, \quad h_{K_{\s}}:= \mbox{diameter of cell}\,\, K_{\s}.
\]
Let $W_{h} \subset W$, a  finite element space of piecewise constant polynomial, given by
\[
W_{h}=\left\{v:v|_{K'}=c_{K_{\s}},\,\,\forall \, K_{\s} \in S^{2}_{h}\right\}.
\]
For the spatial discretization, let $\Omega_{h}$ be a family of shape regular triangulation  of
the domain $\Omega$. 
Further, the mesh size  is denoted   by
\[
h_{\x}:= \max_{K \in \Omega_{h}} h_{\x,K},\qquad  h_{\x,K}:= \mbox{diameter of the cell}\,\, K.
\]
And the finite
element space    of piecewise linear polynomials  $V_{h} $ that vanish on the inlet boundary $  \partial \Omega_{\s,-} $ is defined as
\[
V_{h}:=\left\{v \in \mathbb{C}(\overline{\Omega}):v_{h_{\x}}|_{K} \in \mathbb{P}_{1}(K),\,\,\forall~ K \in \Omega_{h},\,\,v_{h_{\x}}|_{\partial \Omega_{\s,-}}=0 \right\}.
\]
For $u,v\in V_h$, the stabilized SUPG bilinear form  is given by
\[
a_{SUPG}(u,v)=a(u,v)+\sum_{K \in \Omega_{h}}\delta_{K}\left(\s \cdot \nabla u, \s \cdot \nabla v \right)_{K}, 
\]
where  $\delta_{K} > 0$ is an user chosen stabilization parameter. For the convergence analysis, we assume that
\begin{equation}\label{stbasm}
    0 \leq \delta_{K} \leq \delta_{0} h,\quad \delta_{0} > 0.
\end{equation}
Further, the corresponding SUPG-norm is given by
\[
\|v_{h}\|^{2}_{SUPG}:=\left( \sum_{K \in \Omega_{h}} \delta_{K} \|\s \cdot \nabla v_{h} \|^{2}_{L^{2}(K)}+\||\s \cdot \mathbf{n}|^{1/2}v_{h}\|^{2}_{L^{2}(\partial \Omega_{s,+})}\right).
\]
Moreover, the bilinear form associated with the SUPG discretization
is coercive with respect to the $\|\cdot\|_{SUPG}$ by means of (\ref{eu4}). 

\subsection{Semi-discrete method}
Let $\{\phi_{i}\}$ and $\{\psi_{l}\}$ be the basis of the finite dimensional spaces $W_h$ and $V_h$, respectively, \emph{i.e},
\[
W_{h}=\mbox{span}\{\phi_{i}\},\,\, i=1,2,\ldots,N_{\s},\qquad V_{h}=\mbox{span}\{\psi_{l}\},\,\, l=1,2,\ldots,N_{\x}.
\]
Then,   the finite element  space  $ \mathcal{P}^{1,0}_{h} $ is defined as   
\[
\mathcal{P}^{1,0}_{h}:=W_{h}\otimes V_{h}=\bigg\{\zeta \colon \zeta= \sum_{l=1}^{N_\s}\sum_{i=1}^{N_\x} \zeta_{i l}\phi_{i}\psi_{l},\,\,\lambda_{i l} \in \mathbb{R} \bigg\}.
\]
Any discrete function  $v_h \in \mathcal{P}^{1,0}_{h}$ is given by
\[
v_h=\sum_{i=1}^{N_\s}\sum_{l=1}^{N_\x} v_{i l}^n\phi_{i}(\s)\psi_{l}(\x)
\]
and associated advection operator $\s \cdot \nabla v_h$ is expressed as
\[
\s \cdot \nabla v_h = \s \cdot \sum_{i=1}^{N_\s} \sum_{l=1}^{N_\x}v_{i l}^n\phi_{i}(\s)\nabla \psi_{l}(\x).
\]
The efficient way to handle these entries in the associated mass and stiffness  matrices from the resulting  finite element approximation is presented in Section \ref{sec4}. 

Now, by
using the finite element  space  $ \mathcal{P}^{1,0}_{h} $, the respective semi-discrete form   \eqref{fe1} and \eqref{fe2} read:

\noindent
\textbf{Step 1.} Find $\tilde{u}_{h_{\x},h_{\s}}:(t^{n},t^{n+1})\rightarrow \mathcal{P}^{1,0}_{h}$ with $\tilde{u}_{h_{\x},h_{\s}}(t^{n})=u_{h_{\x},h_{\s}}(t^{n})$ such that
\begin{equation}\label{sfe1}
\int_{S^{2}}(\tilde{u}_{t,h_{\x},h_{\s}},\zeta)_{\x}+\int_{S^{2}}(\mathcal{K}_{\sigma}\tilde{u}_{h_{\x},h_{\s}},\zeta)_{\x}=0,\,\,\forall \zeta \in \mathcal{P}^{1,0}_{h}.
\end{equation}
\textbf{Step 2.} Find $u_{h_{\x},h_{\s}}:(t^{n},t^{n+1})\rightarrow \mathcal{P}^{1,0}_{h}$ with $u_{h_{\x},h_{\s}}(t^{n})=\tilde{u}_{h_{\x},h_{\s}}(t^{n+1})$ such that
\begin{equation}\label{sfe2}
\begin{array}{ll}
\displaystyle \int_{S^{2}}(u_{t,h_{\x},h_{\s}},\zeta)_{\x}+\int_{S^{2}}a_{SUPG}(u_{h_{\x},h_{\s}},\zeta)=\displaystyle\int_{S^{2}}(f,\zeta)_{\x}\\[8pt]
\hspace{4cm}\displaystyle+\int_{S^{2}}\sum_{K \in \Omega_{h}}\delta_{K}\left(f-u_{t,h_{\x},h_{\s}}, \s \cdot \nabla \zeta \right)_{K},\,\,\forall \zeta \in \mathcal{P}^{1,0}_{h}.
\end{array}
\end{equation}
To simplify notations, we denote $v_{h_{\x},h_{\s}}$ by $v_h$ and also use similar notations throughout this paper.

\subsection{Temporal discretization}
We consider a uniform partition of the time interval $[0, T]$ with $\Delta t= T/N$, i.e., $t_{n} =n  \Delta t ,\,\, n = 0,1,\ldots, N$. 
After discretizing the temporal variable  by   the implicit backward Euler scheme, the  fully discrete operator-split form of the model problem \eqref{RT1} reads:
\vskip 2mm
\noindent
\textbf{Step 1.} For a given $u^{n}_{h} \in \mathcal{P}^{1,0}_{h}$, find $\tilde{u}^{n+1}_{h} \in  \mathcal{P}^{1,0}_{h}$  such that
\begin{equation}\label{tfe1}
\int_{S^{2}}\left(\partial_{\Delta t} \tilde{u}^{n+1}_{h},\zeta \right)_{\x}+\int_{S^{2}}(\mathcal{K}_{\sigma}\tilde{u}^{n+1}_{h},\zeta)_{\x}=0,\,\,\zeta \in \mathcal{P}^{1,0}_{h},
\end{equation}
where $\partial_{\Delta t} \tilde{u}^{n+1}_{h}=(\tilde{u}^{n+1}_{h}-u^{n}_{h})/\Delta t$.

\noindent
\textbf{Step 2.} Update the solution $\tilde{u}^{n+1}_{h}$ from (\ref{tfe1}) by finding $u^{n+1}_{h}\in \mathcal{P}^{1,0}_{h}$  such that
\begin{equation}\label{tfe2}
\begin{array}{ll}
\displaystyle \int_{S^{2}}\left(\partial_{\Delta t}u^{n+1}_{h},\zeta \right)_{\x}+\int_{S^{2}}a_{SUPG}(u^{n+1}_{h},\zeta)_{\x}=\displaystyle\int_{S^{2}}(f^{n+1},\zeta)_{\x}\\[8pt]
\hspace{2cm}\displaystyle+\int_{S^{2}}\sum_{K \in \Omega_{h}}\delta_{K}\left(f^{n+1}-\partial_{\Delta t}u^{n+1}_{h}, \s \cdot \nabla \zeta \right)_{K}, \,\,\zeta \in \mathcal{P}^{1,0}_{h},
\end{array}
\end{equation}
where $\partial_{\Delta t}u^{n+1}_{h}=(u^{n+1}_{h}-\tilde{u}^{n+1}_{h})/\Delta t$.

\section{A \emph{priori} error estimate : stability and convergence analysis}\label{sec3}
We now discuss the stability and the convergence analysis for the proposed numerical scheme.  We first establish interpolation error estimates and then discuss the local truncation error of the two-step method. After that,  both the local errors are combined  to obtain a global error estimate.

\subsection{Stability result}
The stability estimate of the two-step operator-splitting  method (\ref{tfe1})-(\ref{tfe2}) is derived here.

\begin{theorem}{\label{stablemma}}   Assume that  the stabilization parameter $\delta_{K}$ satisfy
\begin{equation}\label{cond1}
	\delta_{K} \leq \dfrac{\Delta t}{4}, \quad 	\delta = \max\{	\delta_{K}\}, \quad \Delta t \leq \dfrac{1}{2}.
	\end{equation}
	Then, the solutions   $\tilde{u}^{n+1}_{h}$ and $u^{n+1}_{h}$ of the two-step algorithm (\ref{tfe1}) and (\ref{tfe2}) satisfy
\begin{equation}\label{finstab}
 \|u^{n}_{h}\|^{2}_{0}+\Delta t\sum_{m=0}^{n-1} \int_{S^{2}}\|u^{m+1}_{h}\|^{2}_{SUPG} \leq e^{2 T}\left( \|u^{0}_{h}\|^{2}_{0}
+2\Delta t(1+4 \delta \Delta t)\sum_{m=0}^{n-1}  \|f^{m+1}\|^{2}_{0}\right).
\end{equation}
\end{theorem}
\noindent
\textbf{ Proof.}
Setting $\zeta=\tilde{u}^{n+1}_{h}$ in (\ref{tfe1}) and using $2(a-b)a=a^{2}-b^{2}+(a-b)^{2}$ with Lemma \ref{prop_k2}, we obtain
\begin{equation}\label{st1}
\begin{array}{ll}
\frac{1}{2}\|\tilde{u}^{n+1}_{h}\|^{2}_{0}+\frac{1}{2}\|\tilde{u}^{n+1}_{h}-u^{n}_{h}\|^{2}_{0}+\Delta t \underline{\sigma}_{a} \|\tilde{u}^{n+1}_{h}\|^{2}_{0}\leq \frac{1}{2}\|u^{n}_{h}\|^{2}_{0}.
\end{array}
\end{equation}
By  neglecting  the positive terms from the left hand side of above equation  to deduce that 
\begin{equation}\label{st2}
\|\tilde{u}^{n+1}_{h}\|^{2}_{0}\leq \|u^{n}_{h}\|^{2}_{0}.
\end{equation}
Next, setting $\zeta=u^{n+1}_{h}$ in (\ref{tfe2}), we get
\begin{equation}\label{st3}
\begin{array}{ll}
\displaystyle \dfrac{1}{2}\|u^{n+1}_{h}\|^{2}_{0}-\dfrac{1}{2}\|\tilde{u}^{n+1}_{h}\|^{2}_{0}+\dfrac{1}{2}\|u^{n+1}_{h}-\tilde{u}^{n+1}_{h}\|^{2}_{0}+\dfrac{\Delta t}{2} \int_{S^{2}}\|u^{n+1}_{h}\|^{2}_{SUPG} \\[8pt]
\hspace{1cm}\leq \displaystyle \Delta t \left| \int_{S^{2}}(f^{n+1},u^{n+1}_{h})_{\x}\right|+\Delta t\left|  \int_{S^{2}}\sum_{K \in \Omega_{h}}\delta_{K}(f^{n+1},\s \cdot \nabla u^{n+1}_{h})_{K}\right|\\[8pt]
\hspace{2cm}+\displaystyle \left|  \int_{S^{2}}\sum_{K \in \Omega_{h}}\delta_{K}(u^{n+1}_{h}-\tilde{u}^{n+1}_{h},\s \cdot \nabla u^{n+1}_{h})_{K}\right|.
\end{array}
\end{equation}
Using C-S inequality and Young inequality,
first two terms of  the right hand side are bounded by
\begin{equation}\label{st4}
\Delta t \left| \int_{S^{2}}(f^{n+1},u^{n+1}_{h})_{\x}\right| \leq \dfrac{\Delta t}{2} \|f^{n+1}\|^{2}_{0}+ \frac{\Delta t}{2} \|u^{n+1}_{h}\|^{2}_{0}.
\end{equation}
\begin{equation}\label{st5}
\begin{array}{ll}
\displaystyle  \Delta t \left| \int_{S^{2}}\sum_{K \in \Omega_{h}}\delta_{K}(f^{n+1},\s \cdot \nabla u^{n+1}_{h})_{K}\right| & \displaystyle \leq  2\Delta t \int_{S^{2}} \sum_{K \in \Omega_{h}}\delta_{K} \|f^{n+1}\|^{2}_{L^{2}(K)}\\[20pt]
& \displaystyle \hspace{1cm}+\frac{\Delta t}{8}\int_{S^{2}}\|u^{n+1}_{h}\|^{2}_{SUPG}.
\end{array}
\end{equation}
Again,  employing C-S and Young inequalities with assumptions (\ref{cond1}) to deduce that
\begin{equation}\label{st6}
\left|  \int_{S^{2}}\sum_{K \in \Omega_{h}}\delta_{K}(u^{n+1}_{h}-\tilde{u}^{n+1}_{h},\s \cdot \nabla u^{n+1}_{h})_{K}\right|   \leq  \frac{1}{2}\|u^{n+1}_{h}-\tilde{u}^{n+1}_{h}\|^{2}_{0}+\frac{\Delta t}{8}\int_{S^{2}}\|u^{n+1}_{h}\|^{2}_{SUPG}.
\end{equation}
Then,  combing (\ref{st4})-(\ref{st6}), we have
\begin{equation}\label{st7}
\begin{array}{ll}
\displaystyle \hspace{-1cm}\|u^{n+1}_{h}\|^{2}_{0}-\|\tilde{u}^{n+1}_{h}\|^{2}_{0}+\frac{\Delta t}{2} \int_{S^{2}}\|u^{n+1}_{h}\|^{2}_{SUPG} \\[12pt]
\displaystyle \leq \Delta t \|u^{n+1}_{h}\|^{2}_{0}+\Delta t  \|f^{n+1}\|^{2}_{0}+ 4\Delta t  \int_{S^{2}} \sum_{K \in \Omega_{h}}\delta_{K}\|f^{n+1}\|^{2}_{L^{2}(K)}.
\end{array}
\end{equation}
This can be  reduced as follows
\begin{equation}\label{st70}
\begin{array}{ll}
\displaystyle \hspace{-1cm}(1-\Delta t)\|u^{n+1}_{h}\|^{2}_{0}+\frac{\Delta t}{2} \int_{S^{2}}\|u^{n+1}_{h}\|^{2}_{SUPG} \\[12pt]
\displaystyle \leq \|\tilde{u}^{n+1}_{h}\|^{2}_{0} +\Delta t  \|f^{n+1}\|^{2}_{0}+ 4\Delta t  \int_{S^{2}} \sum_{K \in \Omega_{h}}\delta_{K}\|f^{n+1}\|^{2}_{L^{2}(K)},
\end{array}
\end{equation}
using  $1/(1-\Delta t) \leq 1+2 \Delta t \leq 2$. 
Employ  (\ref{st70}) in (\ref{st7}), we deduce that
\begin{equation}\label{st71}
\begin{array}{ll}
\displaystyle 
\|u^{n+1}_{h}\|^{2}_{0}+\Delta t \int_{S^{2}}\|u^{n+1}_{h}\|^{2}_{SUPG} \\[12pt]
\displaystyle \qquad \leq \left(1+2\Delta t \right)\|\tilde{u}^{n+1}_{h}\|^{2}_{0}
+2\Delta t (1+4 \delta \Delta t)\|f^{n+1}\|^{2}_{0}.
\end{array}
\end{equation}
Adding stability results of both the steps  (\ref{st2}) and (\ref{st71}), we have
\begin{equation}\label{st8}
\begin{array}{ll}
\displaystyle \|u^{n+1}_{h}\|^{2}_{0}+\Delta t \int_{S^{2}}\|u^{n+1}_{h}\|^{2}_{SUPG} \\[12pt]
\displaystyle \qquad \leq \left(1+2\Delta t \right)\|u^{n}_{h}\|^{2}_{0}
+2\Delta t (1+4 \delta \Delta t)\|f^{n+1}\|^{2}_{0}.
\end{array}
\end{equation}
Now summing over $m =0,1,\ldots,n-1$, we get that 
\begin{equation}\label{st80}
\begin{array}{ll}
\displaystyle \|u^{n}_{h}\|^{2}_{0}+\Delta t \sum_{m=0}^{n-1} \int_{S^{2}}\|u^{m+1}_{h}\|^{2}_{SUPG} \\[12pt]
\displaystyle \qquad\leq 2\Delta t\sum_{m=0}^{n-1}\|u^{m}_{h}\|^{2}_{0}+ \|u^{0}_{h}\|^{2}_{0}+2\Delta t (1+4 \delta \Delta t)\sum_{m=0}^{n-1}\|f^{m+1}\|^{2}_{0}.
\end{array}
\end{equation}
By using  discrete Grownwall's lemma, we  obtain the stated stability result of the lemma. 
\eop \hfill

\noindent\remark
In Theorem \ref{stablemma}, the stability condition of the discrete method (\ref{tfe1}) and (\ref{tfe2}) is established with the stability parameter $\delta_{K}$ satisfies $\delta_{K}=\mathcal{O}(\Delta t)$. From (\ref{stbasm}), we would be able to take $\Delta t \sim h$. 
For more details on the choice of stabilization parameter,  one may see the detailed discussion in \cite{john2008finite,Eril_CMAME1}. 

\subsection{Convergence analysis}
In this subsection,  error approximation for  the numerical solution.
To derive the error estimate of the operator-splitting finite element discretization(\ref{tfe1})-(\ref{tfe2}), we denote 
\[
\displaystyle \Pi_{h}=\pi_{h_{\x}}\pi_{h_{\s}}=\pi_{h_{\s}}\pi_{h_{\x}}.
\]
Here, $ \pi_{h_{\x}}v \in V_{h}$, the elliptic projection of $v \in V$ and $ \pi_{h_{\s}}w \in W_{h}$, the angular interpolant  of $w \in W$. 
By applying the argumentation form \cite[Lemma  4.2]{apel2005clement}, we have 
\begin{equation}\label{sinterperr}
\|w-\pi_{h_{\s}}w\|^{2}_{L^{2}(S^{2})} \leq C \,\, h^{2}_{\s} |w|^{2}_{H^{1}(S^{2})},\,\,\forall \, w \in  H^{1}(S^{2}).
\end{equation}
Using Galerkin orthogonality,  $ \pi_{h_{\x}}u $ satisfies
\[
a_{SUPG}(\pi_{h_{\x}}u, v_{h_{\x}})=a_{SUPG}(u, v_{h_{\x}}) \,\,\forall \, v_{h_{\x}} \in V_{h}.
\]
Applying as in \cite{johnson1984finite},  we have
\begin{equation}\label{xinterperr}
\|u-\pi_{h_{\x}}u\|^{2}_{SUPG} \leq C\, h_{\x}^{3}\|u\|^{2}_{H^{2}(\Omega)},\,\,\forall \,u \in V \cap H^{2}(\Omega).
\end{equation}

\vskip 2mm Now, we discuss the error estimate of the discrete problems (\ref{tfe1})-(\ref{tfe2}).
The local truncation error in  the first step (\ref{tfe1}) is denoted by 
\[
\tilde{E}^{n}_{h}=\tilde{u}^{n}_{h}-\Pi_{h}\tilde{u}(t^{n}),
\]
where $ \tilde{u}^{n}_{h} $ is the fully discrete solution of  the first step (\ref{tfe1}) and $ \tilde{u} $ is  the weak solution of (\ref{fe1}). The error term $\tilde{E}^{n+1}_{h}$ solves the following equation
\begin{equation}\label{ce1}
\int_{S^{2}}(\partial_{\Delta t} \tilde{E}^{n+1}_{h},\zeta)_{\x}+\int_{S^{2}}\mathcal{A}(\tilde{E}^{n+1}_{h},\zeta)=\int_{S^{2}}(I_{1},\zeta)_{\x}+\int_{S^{2}}\mathcal{A}(I_{2},\zeta),\,\, \zeta \in \mathcal{P}^{1,0}_{h},
\end{equation}
where 
\[
I_{1}=\tilde{u}_{t}(t^{n+1})-\partial_{\Delta t}\Pi_{h}\tilde{u}(t^{n+1}),\quad I_{2}=\tilde{u}(t^{n+1})-\Pi_{h}\tilde{u}(t^{n+1}),
\]
and $\mathcal{A}(v,w) =(\mathcal{K}_{\sigma}v,w)_{\x},$ 
for any discrete function $v$ and $w$. 

\vskip 2mm
Next, we  discuss the error estimates  
for both  steps (\ref{tfe1}) and (\ref{tfe2}) subsequently in the upcoming lemmas.
\begin{lemma}\label{step1err1}
	The local truncation error $ \tilde{E}^{n+1}_{h} $ associated with the angular discretization satisfies
	\begin{equation}\label{cestep1}
 \sum_{m=0}^{n-1}\left(\|\tilde{E}^{m+1}_{h}\|^{2}_{0}-\|E^{m}_{h}\|^{2}_{0}\right)
	\leq C \Delta t \bigg[\Delta t \int_{0}^{T}\|\tilde{u}_{tt}\|^{2}_{0}+ h_{\s}^{2}\|\tilde{u}_{t}\|^{2}_{H^{1}(L^{2})}  + h_{\x}^{3}\|\tilde{u}_{t}\|^{2}_{L^{2}(H^{2})} \bigg].
	\end{equation}	
\end{lemma}
\noindent
\textbf{Proof.}
Setting $\zeta=\tilde{E}^{n+1}_{h}$ in (\ref{ce1}), we get
\begin{equation}\label{ce20}
\begin{array}{ll}
\displaystyle \hspace{-2cm}
\frac{1}{2}\|\tilde{E}^{n+1}_{h}\|^{2}_{0}-\frac{1}{2}\|E^{n}_{h}\|^{2}_{0}+\Delta t \int_{S^{2}}\mathcal{A}\left(\tilde{E}^{n+1}_{h},\tilde{E}^{n+1}_{h}\right) \\[12pt]
\displaystyle \leq \Delta t  \int_{S^{2}}\left|\left(I_{1},\tilde{E}^{n+1}_{h}\right)_{\x}\right|+\Delta t \int_{S^{2}}\left|\mathcal{A}\left(I_{2},\tilde{E}^{n+1}_{h}\right)\right|.
\end{array}
\end{equation}
By using the argumentation  from  (\ref{st1}), it can be deduced that
\begin{equation}\label{ce2}
\frac{1}{2}\|\tilde{E}^{n+1}_{h}\|^{2}_{0}-\frac{1}{2}\|E^{n}_{h}\|^{2}_{0}+\Delta t  \underline{\sigma}_{a} \|\tilde{E}^{n+1}_{h}\|^{2}_{0}\leq \Delta t  \int_{S^{2}}\left|\left(I_{1},\tilde{E}^{n+1}_{h}\right)_{\x}\right|+\Delta t \int_{S^{2}}\left|\mathcal{A}\left(I_{2},\tilde{E}^{n+1}_{h}\right)\right|.
\end{equation}
Let us consider the  first term on  right hand side. Applying C-S and Young inequalities, we deduce that 
\begin{equation}\label{ce3}
\begin{array}{ll}
\displaystyle \Delta t \int_{S^{2}}\left|\left(I_{1},\tilde{E}^{n+1}_{h}\right)_{\x}\right| & \leq 4 \Delta t \|\tilde{u}_{t}(t^{n+1})-\partial_{\Delta t}\Pi_{h}\tilde{u}(t^{n+1})\|^{2}_{0}+\dfrac{\Delta t}{16}\|\tilde{E}^{n+1}_{h}\|^{2}_{0}\\[8pt]
& \leq 8\Delta t\|\Pi_{h}\tilde{u}_{t}(t^{n+1})-\partial_{\Delta t}\Pi_{h}\tilde{u}(t^{n+1})\|^{2}_{0}\\[8pt]
&
\hspace{1cm}+8\Delta t\|\tilde{u}_{t}(t^{n+1})-\Pi_{h}\tilde{u}_{t}(t^{n+1})\|^{2}_{0}+\dfrac{\Delta t}{16}\|\tilde{E}^{n+1}_{h}\|^{2}_{0}\\[8pt]
& \displaystyle \leq C \Delta t^{2} \int_{t^{n} }^{t^{n+1}}\|\Pi_{h}\tilde{u}_{t t}(t^{n+1})\|^{2}_{0}+C\Delta t \|\tilde{u}_{t}(t^{n+1})-\pi_{h_{\s}}\tilde{u}_{t}(t^{n+1})\|^{2}_{0}\\[8pt]
& \displaystyle \hspace{1cm} +C\Delta t\|\pi_{h_{\s}}\tilde{u}_{t}(t^{n+1})-\Pi_{h}\tilde{u}(t^{n+1})\|^{2}_{0}+\dfrac{\Delta t}{16}\|\tilde{E}^{n+1}_{h}\|^{2}_{0}\\[8pt]
& \displaystyle \leq C \Delta t^{2} \int_{t^{n} }^{t^{n+1}}\|\Pi_{h}\tilde{u}_{t t}(t^{n+1})\|^{2}_{0}+C \Delta t \,h_{\s}^{2}\|\tilde{u}_{t}\|^{2}_{H^{1}(L^{2})} \\[8pt]
&  \hspace{1cm}
+C \Delta t \, h_{\x}^{3}\|\tilde{u}_{t}\|^{2}_{L^{2}(H^{2})}+\dfrac{\Delta t}{16}\|\tilde{E}^{n+1}_{h}\|^{2}_{0}.
\end{array}
\end{equation}
Next, the second term is decomposed as 
\[
\begin{array}{ll}
\displaystyle \int_{S^{2}}\mathcal{A}\left(I_{2},\tilde{E}^{n+1}_{h}\right)=\int_{S^{2}}\mathcal{A}\left(\tilde{u}(t^{n+1})-\pi_{h_{\s}}\tilde{u}(t^{n+1}),\tilde{E}^{n+1}_{h}\right)\\[8pt] 
\hspace{4cm}\displaystyle
+\int_{S^{2}}\mathcal{A}\left(\pi_{h_{\s}}\tilde{u}(t^{n+1})-\Pi_{h}\tilde{u}(t^{n+1}),\tilde{E}^{n+1}_{h}\right).
\end{array}
\]
Using Lemma \ref{prop_k2} with the interpolation and projection estimates, the first term can be deduced by means of C-S and Young inequalities,
\begin{equation}\label{ce4}
\begin{array}{ll}
\displaystyle \Delta t \int_{S^{2}}\left|\mathcal{A}\left(\tilde{u}(t^{n+1})-\pi_{h_{\s}}\tilde{u}(t^{n+1}),\tilde{E}^{n+1}_{h}\right)\right| &\\[8pt] 
\hspace{2cm}\displaystyle \leq \Delta t (2\overline{\sigma}_{\s}+\overline{\sigma}_{a}) \left\|\tilde{u}(t^{n+1})-\pi_{h_{\s}}\tilde{u}(t^{n+1})\right\|_{0} \|\tilde{E}^{n+1}_{h}\|_{0}\\[8pt]
\hspace{2cm} \leq C \Delta t\, h_{\s}^{2}\|\tilde{u}\|^{2}_{H^{1}(L^{2})}+\dfrac{\Delta t}{32}\|\tilde{E}^{n+1}_{h}\|^{2}_{0}.
\end{array}
\end{equation}
Similarly, one can obtain that
\begin{equation}\label{ce5}
\Delta t \int_{S^{2}}\left|\mathcal{A}\left(\pi_{h_{\s}}\tilde{u}\left(t^{n+1}\right)-\Pi_{h}\tilde{u}\left(t^{n+1}\right),\tilde{E}^{n+1}_{h}\right)\right| \leq C \Delta t\, h_{\x}^{3}\|\tilde{u}\|^{2}_{L^{2}(H^{2})}+\dfrac{\Delta t}{32}\|\tilde{E}^{n+1}_{h}\|^{2}_{0}.
\end{equation}
Combing estimates (\ref{ce4}) and (\ref{ce5}), we get 
\begin{equation}\label{ce51}
\Delta t \int_{S^{2}}\left|\mathcal{A}\left(I_{2},\tilde{E}^{n+1}_{h}\right)\right| \leq C \Delta t\, h_{\s}^{2}\|\tilde{u}\|^{2}_{H^{1}(L^{2})}+ C \Delta t\,h_{\x}^{3}\|\tilde{u}\|^{2}_{L^{2}(H^{2})}+\dfrac{\Delta t}{16}\|\tilde{E}^{n+1}_{h}\|^{2}_{0}.
\end{equation}
Next, employing (\ref{ce3}) and (\ref{ce51}) in (\ref{ce2}), it can be devised that 
\begin{equation}\label{ce52}
\begin{array}{ll}
\displaystyle 
\frac{1}{2}\|\tilde{E}^{n+1}_{h}\|^{2}_{0}
-\frac{1}{2}\|E^{n}_{h}\|^{2}_{0} +  \Delta t \underline{\sigma}_{a} \|\tilde{E}^{n+1}_{h}\|^{2}_{0} \leq C \Delta t^{2} \int_{t^{n} }^{t^{n+1}}\|\Pi_{h}\tilde{u}_{t t}(t^{n+1})\|^{2}_{0}\\[12pt]
\displaystyle \hspace{2cm} +C \Delta t \,h_{\s}^{2}\|\tilde{u}_{t}\|^{2}_{H^{1}(L^{2})}+C \Delta t \, h_{\x}^{3}\|\tilde{u}_{t}\|^{2}_{L^{2}(H^{2})}+\dfrac{\Delta t}{8}\|\tilde{E}^{n+1}_{h}\|^{2}_{0}.
\end{array}
\end{equation}
By using assumption \ref{Asm2}, it is reduced as 
\begin{equation}\label{ce53}
\begin{array}{ll}
\displaystyle 
\frac{1}{2}\|\tilde{E}^{n+1}_{h}\|^{2}_{0}
-\frac{1}{2}\|E^{n}_{h}\|^{2}_{0}  \leq C \Delta t^{2} \int_{t^{n} }^{t^{n+1}}\|\Pi_{h}\tilde{u}_{t t}(t^{n+1})\|^{2}_{0}\\[12pt]
\displaystyle \hspace{2cm} +C \Delta t \,h_{\s}^{2}\|\tilde{u}_{t}\|^{2}_{H^{1}(L^{2})}+C \Delta t \, h_{\x}^{3}\|\tilde{u}_{t}\|^{2}_{L^{2}(H^{2})}.
\end{array}
\end{equation}
Finally,  summing over $m=1,2,\ldots,n-1$, we obtain the required results.
\hfill \eop

\noindent\remark
Note that  the coefficient of $\|\tilde{E}^{n+1}_{h}\|^{2}_{0}$ in the left side of (\ref{ce52}) can be taken smaller as per our convenience  by the means of C-S and Young inequalities. 


\vskip 2mm
Next, we discuss the bound for the local truncation error term for the second step,  which is given by
\[
E^{n}_{h}= u^{n}_{h}-\Pi_{h}u(t^{n}).
\]
Let $\pi_{h_{\x}}u \in V_{h}$ be the elliptic projection of $u \in V$, we have
\[
a_{SUPG}(\pi_{h_{\x}}u,v_h)=a_{SUPG}(u,v_h), \,\, \forall \,v_{h} \in V_{h}.
\]
Then, one can claim that
\[
a_{SUPG}(\pi_{h_{\s}}\pi_{h_{\x}}u,v_h)=\left(\pi_{h_{\s}}f-\pi_{h_{\s}}u_{t},v_{h}+\delta \s \cdot \nabla v_{h}\right), \,\, \forall \, v_{h} \in V_{h}.
\]
And, we have 
\begin{equation}\label{ce6}
\int_{S^{2}}a_{SUPG}(\Pi_{h}u,\zeta)=\int_{S^{2}}\left(\pi_{h_{\s}}f-\pi_{h_{\s}}u_{t},\zeta+\delta \s \cdot \nabla \zeta\right)_{\x},\,\,\forall \, \zeta \in \mathcal{P}^{1,0}_{h}.
\end{equation}
The error term $ E^{n+1}_{h} $ satisfies the following equation
\begin{equation}\label{ce7}
\begin{array}{ll}
\displaystyle  \int_{S^{2}}(\partial_{\Delta t} E^{n+1}_{h},\zeta)_{\x}+\int_{S^{2}}a_{SUPG}(E^{n+1}_{h},\zeta) \\ [8pt]
\displaystyle
= \int_{S^{2}}\left(\Lambda^{n+1}_{u}+\Lambda^{n+1}_{f},\zeta+\delta \s \cdot \nabla \zeta\right)_{\x}  - \delta \int_{S^{2}}(\partial_{\Delta t} E^{n+1}_{h},\s \cdot \nabla \zeta)_{\x} ,\,\,\forall \zeta \in \mathcal{P}^{1,0}_{h},
\end{array}
\end{equation}
where functions $ \Lambda^{n+1}_{u} $ and $ \Lambda^{n+1}_{f} $ are defined as 
\[
\Lambda^{n+1}_{u}=\left(\pi_{h_{\s}}u_{t}(t^{n+1})-\partial_{\Delta t}u(t^{n+1})
\right)\quad \mbox{and} \quad \Lambda^{n+1}_{f}=f^{n+1}-\pi_{h_{\s}}f(t).
\]

\noindent
\begin{lemma}\label{step2err1}
	The local truncation error $ E^{n+1}_{h} $ associated with the spatial discretization satisfies
	\begin{equation}\label{cestep2}
	\begin{array}{ll}
	\displaystyle
	\sum_{m=0}^{n-1}\left(\|E^{m+1}_{h}\|^{2}_{0}-\|\tilde{E}^{m+1}_{h}\|^{2}_{0}\right)+\dfrac{\Delta t}{2}\sum_{m=0}^{n-1}\int_{S^{2}}\|E^{m+1}_{h}\|_{SUPG}^{2}\\ [14pt]
	\displaystyle \hspace{1cm} \leq C \bigg[\Delta t^{2} \int_{0}^{T}\left(\|u_{tt}\|^{2}_{0}+\|u_{ttt}\|^{2}_{0}\right) \\ [14pt]
	\displaystyle \hspace{2cm}
	+ \Delta t \, h_{\x}^{3}\sum_{m=0}^{n-1}\left(\|u\|^{2}_{L^{2}(H^{2})}+\|u_{t}\|^{2}_{L^{2}(H^{2})}+\|u_{t  t}\|^{2}_{L^{2}(H^{2})}\right)\\ [14pt]
	\displaystyle \hspace{2cm}+\Delta t \, h_{\s}^{2}\sum_{m=0}^{n-1}\left(\|u\|^{2}_{H^{1}(H^{1})}+\|u_{t}\|^{2}_{H^{1}(H^{1})}+\|u_{t t}\|^{2}_{H^{1}(H^{1})}\right) \bigg].
	\end{array}
	\end{equation}
\end{lemma}
\noindent
\textbf{Proof.}
Setting $\zeta = E^{n+1}_{h}$ in (\ref{ce7}), we obtain 
\[
\begin{array}{ll}
\displaystyle \hspace{-5mm} \frac{1}{2 \Delta t}\left(\|E^{n+1}_{h}\|^{2}_{0}-\|\tilde{E}^{n+1}_{h}\|^{2}_{0}\right)+\dfrac{1}{2}\int_{S^{2}}\|E^{n+1}_{h}\|_{SUPG}^{2} \\ [12pt]
\displaystyle  \leq \int_{S^{2}} \left|\left(\Lambda^{n+1}_{u}+\Lambda^{n+1}_{f},E^{n+1}_{h}+\delta \s \cdot \nabla E^{n+1}_{h} \right)_{\x}\right| + \delta \int_{S^{2}}\left|(\partial_{\Delta t} E^{n+1}_{h},\s \cdot \nabla E^{n+1}_{h})_{\x} \right|.
\end{array}
\]
Applying C-S inequality,   Young inequality and summing over $m=0,1,\ldots,n-1$, we get
\begin{equation}\label{ce8}
\begin{array}{ll}
\displaystyle
\sum_{m=0}^{n-1}\left(\|E^{m+1}_{h}\|^{2}_{0}-\|\tilde{E}^{m+1}_{h}\|^{2}_{0}\right)+\dfrac{\Delta t}{2}\sum_{m=0}^{n-1}\int_{S^{2}}\|E^{m+1}_{h}\|_{SUPG}^{2}\\ [12pt]
\displaystyle \hspace{2cm} \leq  10(1+\delta)\sum_{m=0}^{n-1}\Delta t \left(\|\Lambda^{m+1}_{u}\|^{2}_{0}+\|\Lambda^{m+1}_{f}\|^{2}_{0}\right)+10 \delta \sum_{m=0}^{n-1}\Delta t \|\partial_{\Delta t}E^{m+1}_{h}\|^{2}_{0}.
\end{array}
\end{equation}
Next, we discuss estimate for last term of the right hand side, \emph{i.e.}, $\|\partial_{\Delta t}E^{m+1}_{h}\|^{2}_{0}$.  
Consider  $\chi^{n+1}=\partial_{\Delta t} E^{n+1}_{h}$ in (\ref{ce7}), it satisfies the following equation
\begin{equation}\label{ce9}
\begin{array}{ll}
\displaystyle
\int_{S^{2}}\left(\partial_{\Delta t} \chi^{n+1},\zeta\right)_{\x}+\int_{S^{2}}a_{SUPG}(\chi^{n+1},\zeta) \\[16pt]
=\displaystyle
\int_{S^{2}}\left(\partial_{\Delta t}\Lambda^{n+1}_{u}+\partial_{\Delta t}\Lambda^{n+1}_{f},\zeta+\delta \s \cdot \nabla \zeta\right)_{\x} - \delta \int_{S^{2}}\left(\partial_{\Delta t} \chi^{n+1},\s \cdot \nabla \zeta\right)_{\x},
\end{array}
\end{equation}
where $\partial_{\Delta t}\chi^{n+1}=(\chi^{n+1}-\chi^{n})/\Delta t$ and $\partial_{\Delta t}\Lambda^{n+1}_{u}$, $\partial_{\Delta t}\Lambda^{n+1}_{f}$ are defined in  a similar way.

Further, assigning  $\zeta=\chi^{n+1}+\delta \partial_{\Delta t} \chi^{n+1}$ in (\ref{ce9}) and after small simplification,  we deduce 
\begin{equation}\label{ce10}
\begin{array}{ll}
\dfrac{1}{2 \Delta t}\bigg(\|\chi^{n+1}\|^{2}_{0}-\|\chi^{n}\|^{2}_{0}\bigg)+\delta \|\partial_{\Delta t}\chi^{n+1}\|^{2}_{0} + \delta \|\s \cdot \nabla \chi^{n+1}\|^{2}_{0}+\dfrac{1}{2}\|(\s \cdot \mathbf{n}) \chi^{n+1}\|^{2}_{L^{2}(\Gamma_{+})}\\[16pt]
\displaystyle +\dfrac{\delta^{2}}{2 \Delta t}\bigg(\|\s \cdot \nabla\chi^{n+1}\|^{2}_{0}-\|\s \cdot \nabla\chi^{n}\|^{2}_{0}\bigg)   +   2\delta \int_{S^{2}}\bigg(\partial_{\Delta t}\chi^{n+1}, \s \cdot \nabla \chi^{n+1}\bigg)_{\x}\\[16pt]
\displaystyle 
\hspace{1cm} \leq   \int_{S^{2}}\bigg(\partial_{\Delta t}\Lambda^{n+1}_{u}+\partial_{\Delta t}\Lambda^{n+1}_{f},\chi^{n+1}+\delta \s \cdot \nabla \chi^{n+1}\bigg)_{\x}\\[16pt]
\displaystyle \hspace{1.5cm}+  \delta \int_{S^{2}}\bigg(\partial_{\Delta t}\Lambda^{n+1}_{u}+\partial_{\Delta t}\Lambda^{n+1}_{f},\partial_{\Delta t} \chi^{n+1}+\delta \s \cdot \nabla \partial_{\Delta t} \chi^{n+1}\bigg)_{\x}
\\[16pt]
\displaystyle \hspace{1.5cm}
-\delta^{2} \int_{S^{2}}\big(\partial_{\Delta t}\chi^{n+1}, \s \cdot \nabla (\partial_{\Delta t}\chi^{n+1})\big)_{\x}.
\end{array}
\end{equation}
Next, we define the  norm $\|\cdot\|_{\x \s}$ as follows:
\[
\|\chi^{n+1}\|_{\x \s}^{2}\equiv  \delta\left\|\s \cdot \nabla\chi^{n+1}+\partial_{\Delta t}\chi^{n+1}\right\|_{0}^{2}+\dfrac{1}{2}\|(\s \cdot \mathbf{n}) \chi^{n+1}\|^{2}_{L^{2}(\Gamma_{+})}.
\] 
By using the above norm definition, we deduce 
\begin{equation}\label{ce11}
\begin{array}{ll}
\displaystyle \dfrac{1}{2 \Delta t}\left(\|\chi^{n+1}\|^{2}_{0}-\|\chi^{n}\|^{2}_{0}\right)+\|\chi^{n+1}\|^{2}_{\x \s}+\dfrac{\delta^{2}}{2 \Delta t}\left(\|\s \cdot \nabla\chi^{n+1}\|^{2}_{0}-\|\s \cdot \nabla\chi^{n}\|^{2}_{0}\right)\\ [12pt]
\displaystyle \hspace{1cm}
\leq  \int_{S^{2}}\left(\partial_{\Delta t}\Lambda^{n+1}_{u}+\partial_{\Delta t}\Lambda^{n+1}_{f},\chi^{n+1}+\delta \s \cdot \nabla \chi^{n+1}\right)_{\x}\\ [12pt]
\displaystyle \hspace{1.2cm} +\delta  \int_{S^{2}}\left(\partial_{\Delta t}\Lambda^{n+1}_{u}+\partial_{\Delta t}\Lambda^{n+1}_{f},\partial_{\Delta t} \chi^{n+1}+\delta \s \cdot \nabla \partial_{\Delta t} \chi^{n+1}\right)_{\x}\\[16pt]
\displaystyle \hspace{1.2cm}
-\delta^{2} \bigg(\partial_{\Delta t}\chi^{n+1}, \s \cdot \nabla (\partial_{\Delta t}\chi^{n+1})\bigg)_{\x}.
\end{array}
\end{equation}
By using C-S and Young inequalities, we have estimated  the right side terms of (\ref{ce11}). After small mathematical simplification, we   obtain the  following inequality
\begin{equation}\label{ce110}
\begin{array}{ll}
\|\chi^{n+1}\|^{2}_{0}-\|\chi^{n}\|^{2}_{0}+\dfrac{\Delta t}{2}\|\chi^{n+1}\|^{2}_{\x \s}+\delta^{2}\left(\|\s \cdot \nabla\chi^{n+1}\|^{2}_{0}-\|\s \cdot \nabla\chi^{n}\|^{2}_{0}\right)\\ [12pt]
\displaystyle \hspace{1cm}
\leq \Delta t\|\chi^{n+1}\|^{2}_{0}+C \Delta t \left(\|\partial_{\Delta t}\Lambda^{m+1}_{u}\|^{2}_{0}+\|\partial_{\Delta t}\Lambda^{m+1}_{f}\|^{2}_{0}\right).
\end{array}
\end{equation}
By summing over $m=1,\ldots,n-1$ and repeating the argument from Theorem \ref{stablemma}, we deduce that
\begin{equation}\label{ce12}
\begin{array}{ll}
\displaystyle 
\|\chi^{n}\|^{2}_{0}+\Delta t \sum_{m=1}^{n-1}\|\chi^{m+1}\|^{2}_{\x \s}+\delta^{2}\|\s \cdot \nabla\chi^{n}\|^{2}_{0} \leq C \exp^{2T} \bigg(\|\chi^{1}\|^{2}_{0} \\[12pt]
\hspace{3cm} \displaystyle  
+C \delta^{2}\|\s \cdot \nabla\chi^{1}\|^{2}_{0} +  C \Delta t \sum_{m=1}^{n-1} \left(\|\partial_{\Delta t}\Lambda^{m+1}_{u}\|^{2}_{0}+\|\partial_{\Delta t}\Lambda^{m+1}_{f}\|^{2}_{0}\right)\bigg).
\end{array}
\end{equation}
Now substituting the  estimate of $\|\partial_{\Delta t}E^{m+1}_{h}\|^{2}_{0}$ from (\ref{ce12}) in (\ref{ce8}) and using the above inequality to get
\begin{equation}\label{ce13}
\begin{array}{ll}
\displaystyle
\sum_{m=0}^{n-1}\left(\|E^{m+1}_{h}\|^{2}_{0}-\|\tilde{E}^{m+1}_{h}\|^{2}_{0}\right)+\dfrac{\Delta t}{2}\sum_{m=0}^{n-1}\int_{S^{2}}\|E^{m+1}_{h}\|_{SUPG}^{2} \leq C \delta \|\partial_{\Delta t} E^{1}_{h}\|^{2}_{0}\\ [12pt]
\displaystyle \hspace{6mm} + C\Delta t \left[\sum_{m=0}^{n-1} \left(\|\Lambda^{m+1}_{u}\|^{2}_{0}+\| \Lambda^{m+1}_{f}\|^{2}_{0}\right)
+\delta \sum_{m=0}^{n-1} \left(\|\partial_{\Delta t}\Lambda^{m+1}_{u}\|^{2}_{0}+\|\partial_{\Delta t} \Lambda^{m+1}_{f}\|^{2}_{0}\right)\right] .
\end{array}
\end{equation}
Applying the standard interpolation results and Taylor's theorem with integral
remainder term, we deduce 
\begin{equation}\label{ce14}
\Delta t \sum_{m=0}^{n-1} \|\Lambda^{m+1}_{u}\|^{2}_{0} \leq C \Delta t^{2} \sum_{m=0}^{n-1} \int_{t^{n} }^{t^{n+1} }\|\Pi_{h}u_{t t}\|^{2}_{0}+C\Delta t \, h_{\x}^{3}\sum_{m=0}^{n-1} \|\pi_{h_{\s}}u_{t}(t^{m+1})\|^{2}_{L^{2}(H^{2})}.
\end{equation}
\begin{equation}\label{ce15}
\begin{array}{ll}
\hspace{-1cm}\displaystyle \Delta t \sum_{m=0}^{n-1} \|\partial_{\Delta t} \Lambda^{m+1}_{u}\|^{2}_{0} \\[8pt]
\displaystyle \leq C \Delta t^{2} \sum_{m=0}^{n-1} \int_{t^{n} }^{t^{n+1} }\|\Pi_{h}u_{ttt}\|^{2}_{0}+C\Delta t \, h_{\x}^{3}\sum_{m=0}^{n-1} \|\pi_{h_{\s}}u_{tt}(t^{m+1})\|^{2}_{L^{2}(H^{2})}.
\end{array}
\end{equation}
In a similar way, we have
\begin{eqnarray}
\Delta t \sum_{m=0}^{n-1} \|\Lambda^{m+1}_{f}\|^{2}_{0} & \leq & C \Delta t \, h_{\s}^{2}\sum_{m=0}^{n-1} \|\partial_{\s} f(t^{m+1})\|^{2}_{0}.\label{ce16} \\[12pt]
\Delta t \sum_{m=0}^{n-1} \|\partial_{\Delta t} \Lambda^{m+1}_{f}\|^{2}_{0} & \leq  & C \Delta t^{2} \int_{0}^{T}\|f_{tt}\|^{2}_{0}+C \Delta t \, h_{\s}^{2}\sum_{m=0}^{n-1} \|\partial_{\s} f_{t}(t^{m+1})\|^{2}_{0}.\label{ce17}
\end{eqnarray}

Next, we need to evaluate bounds for term $\|\partial_{\Delta t} E^{1}_{h}\|^{2}_{0}$. By  using the  argument  from \cite{ahamd_ANM} and following  similar technique as in (\ref{ce12}), we devise that 
\begin{equation}\label{ce140}
\delta \|\partial_{\Delta t} E^{1}_{h}\|^{2}_{0}  \leq C (\Delta t^{2}+h^{3}_{\x}+h^{2}_{\s}).
\end{equation}

By combing the estimates (\ref{ce14})-(\ref{ce140}) in (\ref{ce13}), we obtain 
the  desired result (\ref{cestep2}).
This completes the proof.
\hfill \eop

\vskip 5mm
Next, we present the  convergence estimate of the operator-splitting finite element method (\ref{tfe1})-(\ref{tfe2}) for the model problem (\ref{RT1}).

\begin{theorem}
	The global error  $e^n =u-u^{n}_{h}$ satisfies 
	\begin{equation}\label{mainerr}
	\|e^n\|^{2}_{0}+\dfrac{\Delta t}{2}\sum_{m=0}^{n}\int_{S^{2}}\|e^{m+1}\|^{2}_{SUPG}\leq C(u,\tilde{u})(\Delta t^{2}+h^{3}_{\x}+h^{2}_{\s}).
	\end{equation}
\end{theorem}
\noindent
\textbf{Proof.} By combining the local error estimates (\ref{cestep1}) and (\ref{cestep2}), we obtain that 
\begin{equation}\label{mainerreq1}
    \begin{array}{ll}
	\displaystyle
	\|E^{n}_{h}\|^{2}_{0}+\dfrac{\Delta t}{2}\sum_{m=0}^{n-1}\int_{S^{2}}\|E^{m+1}_{h}\|_{SUPG}^{2}\\ [14pt]
	\displaystyle \hspace{1cm} \leq \|E^{0}_{h}\|^{2}_{0}+C \Delta t \bigg[\Delta t \int_{0}^{T}\|\tilde{u}_{tt}\|^{2}_{0}+ h_{\s}^{2}\|\tilde{u}_{t}\|^{2}_{H^{1}(L^{2})}  + h_{\x}^{3}\|\tilde{u}_{t}\|^{2}_{L^{2}(H^{2})} \bigg]\\ [14pt]
	\displaystyle \hspace{2cm}+C \bigg[\Delta t^{2} \int_{0}^{T}\left(\|u_{tt}\|^{2}_{0}+\|u_{ttt}\|^{2}_{0}\right) \\ [14pt]
	\displaystyle \hspace{2cm}
	+ \Delta t \, h_{\x}^{3}\sum_{m=0}^{n-1}\left(\|u\|^{2}_{L^{2}(H^{2})}+\|u_{t}\|^{2}_{L^{2}(H^{2})}+\|u_{t  t}\|^{2}_{L^{2}(H^{2})}\right)\\ [14pt]
	\displaystyle \hspace{2cm}+\Delta t \, h_{\s}^{2}\sum_{m=0}^{n-1}\left(\|u\|^{2}_{H^{1}(H^{1})}+\|u_{t}\|^{2}_{H^{1}(H^{1})}+\|u_{t t}\|^{2}_{H^{1}(H^{1})}\right) \bigg].
	\end{array}
\end{equation}
Noting that $E^{0}_{h}=0$.  Then,  the estimate  (\ref{mainerreq1}) can be further simplified as
\begin{equation}\label{mainerreq2}
    \begin{array}{ll}
	\displaystyle
	\|E^{n}_{h}\|^{2}_{0}+\dfrac{\Delta t}{2}\sum_{m=0}^{n-1}\int_{S^{2}}\|E^{m+1}_{h}\|_{SUPG}^2 \leq C(u,\tilde{u})(\Delta t^{2}+h^{3}_{\x}+h^{2}_{\s}),
	\end{array}
\end{equation}
where $C(u,\tilde{u})$ is the positive constant, depending upon  $u,\tilde{u}$ in (\ref{mainerreq1}). Finally, by employing the approximation results (\ref{sinterperr}) and (\ref{xinterperr}) and using above argumentation, the  main convergence  result (\ref{mainerr}) is devised. This completes the  proof.
\hfill \eop

\section{Computational Results}\label{sec4}
In this section, we present the numerical algorithm for the proposed discrete scheme and a validation of the theoretical estimates.   

\bigskip
\subsection{Numerical implementation}
We briefly discuss the operator-splitting algorithm for the radiative transfer equation solution, whereas more details on implementation for general scalar equations can be found in~\cite{ganeshan_AMC1} and \cite{barrenechea_ANM2017}. All the numerical experiments are performed in our in-house finite element package~\cite{ganesan2016object,wilbrandt2017parmoon}.

In Section \ref{sec2}, the finite element space $\mathcal{P}^{1,0}_{h}$ is given as follows
\[
\mathcal{P}^{1,0}_{h}:=W_{h}\otimes V_{h}=\bigg\{\lambda \colon \lambda= \sum_{l=1}^{N_s}\sum_{i=1}^{N_x} \lambda_{i l}\phi_{i}\psi_{l},\,\,\lambda_{i l} \in \mathbb{R} \bigg\}.
\]
Further, the discrete solution $u^{n}_{h}\in \mathcal{P}^{1,0}_{h}$ and its gradient are expressed  by
\[
u^{n}_{h}(\x,\s)=\sum_{i=1}^{N_s}\sum_{l=1}^{N_x} u_{i l}^n\phi_{i}(\s)\psi_{l}(\x),\quad \nabla u^{n}_{h}=\sum_{i=1}^{N_s} \sum_{l=1}^{N_x}u_{i l}^n\phi_{i}(\s)\nabla \psi_{l}(\x),
\]
where $ u_{i l}$ are the unknown degrees of freedoms (DOFs).
Define the mass  matrices $M^{1}_{\s},\,M^{2}_{\s} \in \mathbb{R}^{N_{\s} \times N_{\s}}$, where the $(i,j)^{t h}$ entries of these matrices are given by 
\[
(M^{1}_{\s})_{ij}=\int_{S^{2}}\phi_{i}\phi_{j}\,d\s, \qquad  (M^{2}_{\s})_{ij}=\int_{S^{2}}\mathcal{K}\phi_{i}\phi_{j}\,d\s.
\]
Further, the $(l,m)^{t h}$ entries of the matrices $M_{\x},\,A_{\x},\,M^{\delta}_{\x},\, A^{\delta}_{\x}$ are given by
\[
\begin{array}{ll}
\displaystyle (M_{\x})_{lm}=\int_{\Omega}\psi_{l}\psi_{m} \,d\x,\qquad \left(M^{\delta}_{\x}\right)_{lm}=\sum_{K 
	\in \Omega_{h}}\delta_{K}(\psi_{l},\s\cdot \nabla \psi_{m})_{K}\,d\x,\\ 
\displaystyle (A_{\x})_{lm}=\int_{\Omega}\s\cdot \nabla\psi_{l}\psi_{m}\,d\x,\qquad \left(A^{\delta}_{\x}\right)_{lm}=\sum_{K 
	\in \Omega_{h}}\delta_{K}(\s\cdot \nabla \psi_{l},\s\cdot \nabla \psi_{m})_{K}\,d\x,
\end{array}
\]
and $m^{th}$ component of the load vectors $F^{n}_{\x}$ and $ F^{\delta,n}_{\x}$ are given by 
\[
\displaystyle \left(F^{n}_{\x}\right)_{m}=\int_{\Omega}f^{n}\psi_{m}\,d\x,\qquad \left(F^{\delta,n}_{\x}\right)_{m}=\sum_{K 
	\in \Omega_{h}}\delta_{K}(f^{n},\s \cdot \nabla \psi_{m})_{K}.
\]

Here the matrices $ M^{\delta}_{\x},\, A^{\delta}_{\x}$ and $ F^{\delta,n}_{\x}$ are due  to the SUPG stabilization terms.
Further, we use $\widetilde{\mathbf{u}}^{n+1}_{k}$ to denote an array of unknown solution coefficients (DOFs) $\{\tilde{u}^{n+1}_{k,i}\},\,\,k=1,2,\ldots,N_{\x},\,i=1,2,\ldots,N_{\s}$.
With these notations, the system matrix of the $\s$-direction step in the time interval $(t^{n},t^{n+1}]$ becomes:\\

\noindent 
For given $\mathbf{u}^{n}_{k}$, solve
\begin{equation}\label{Ne1}
\mathbf{M}_{\s}\,\widetilde{\mathbf{u}}^{n+1}_{k}=M^{1}_{\s}\mathbf{u}^{n}_{k}, \qquad \mathbf{M}_{\s}=\left(M^{1}_{\s}+ \Delta t \,\sigma_{\tau}M^{1}_{\s}- \Delta t \,\sigma_{s} M^{2}_{\s}\right),
\end{equation}
for $k=1,2,\ldots,N_{\x}$. In this way, we solve numerical solution in the $\s$-direction.

\bigskip

We next discuss the  $\x$-direction
step, where the updated  solution from the $\s$-direction  is used to compute the solution of (\ref{tfe1}). In the $\x$-direction step, we first need to transpose the vector $\widetilde{\mathbf{u}}^{n+1}_{k}$ to obtain $\widetilde{\mathbf{u}}^{n+1}_{\ell}$ and then solve the linear system 
\begin{equation}\label{Ne3}
\left(M_{\x}+M^{\delta}_{\x}+ \Delta t ( A_{\x}+A^{\delta}_{\x})\right)\mathbf{u}^{n+1}_{\ell}  = \Delta t ( F^{n}_{\x}+F^{\delta,n}_{\x})  +(M_{\x}+M^{\delta}_{\x})\widetilde{\mathbf{u}}^{n+1}_{\ell},
\end{equation}
for $\ell=1,2,\ldots,N_{\s}$.  Though the mass matrix $M_\x$ is independent of $\s$, all other matrices in~\eqref{Ne3} depend on $\s$. Therefore, all these $\s$-dependent matrices need to be assembled for each $\ell$ in all time steps. However, the matrix assembling can be avoided for every $\ell$  by assembling and storing the $\s$-dependent matrices in a component form. For example,  $A_{\x}$ can be split as 
\begin{align*}
\displaystyle (A_{\x})_{lm} &= s_1\int_{\Omega_h} \frac{\partial \psi_{l}}{\partial x}\psi_{m}\,d\x + s_2\int_{\Omega_h} \frac{\partial \psi_{l}}{\partial y}\psi_{m}\,d\x +s_3\int_{\Omega_h} \frac{\partial \psi_{l}}{\partial z}\psi_{m}\,d\x\\
& = s_1(A_{\x}^{I})_{lm} + s_2(A_{\x}^{II})_{lm} + s_3(A_{\x}^{III})_{lm}.
\end{align*}

Hence, it is enough to assemble matrices $A_{\x}^{I}$, $A_{\x}^{II}$ and $A_{\x}^{III}$ only once and then multiply it with $s_1$, $s_2$ and $s_3$, respectively, for every $\ell$ in each time step. Following a similar technique for  $\s$-dependent matrices,  matrix assembling for each $\ell$ in every time step can completely be avoided, and it is enough to assemble all these component matrices only once at the beginning of the computation. In this way, we can solve the linear system (\ref{Ne3}) very efficiently.

\bigskip 
\subsection{Validation}
To validate the theoretical estimates discussed in the previous sections, we consider multiple test examples using manufactured solutions.  The scattering phase function $\Phi$ is taken from the previous studies \cite{badri1,wang}. For time discretization, the backward Euler time-stepping method is applied with a final time $T = 1$ and time step $\Delta t=h_{\x}$ in all the numerical experiments.  Furthermore,  the value of  absorption and scattering coefficients are taken as $\sigma_{\tau}=2,\,\,\sigma_{s}=1/2$ for all the test problems.  For validation purposes, we discuss the error estimate of the numerical approximation in  the spatial domain $\mathcal{G}$ and time-domain $(0,1)$ by using the following norm:
\[
\ell_{2}(0,1;L^{2}(\mathcal{G}))=\left(\Delta t\sum_{n=1}^{N}\|u(t^{n})-u^{n}_{h}\|^{2}_{L^{2}(S^{2}_{h}\times \Omega_{h})}\right)^{1/2}
\]
where $L^2$ error is calculated in $\x$- and $\s$-directions, \emph{i.e.}
\[
\|u(t^{n})-u^{n}_{h}\|^{2}_{L^{2}(S^{2}_{h}\times \Omega_{h})}=\int_{S^{2}_{h}}\int_{\Omega_{h}}(u(t^{n})-u^{n}_{h})^{2}d\x d\s.
\]
For all the numerical experiments, the mesh  mesh parameters are taken as $N_{\x}=27, 125, 729, 4913$ and $N_{\s}=48, 192, 768, 3072$.

\begin{example}\label{ex1}
	Consider the model problem~\eqref{RT1} with the exact solution as
	\[
	u(\x,\s,t)=e^{-\alpha t} \sin(\pi x_1)\sin(\pi x_2)\sin(\pi x_3),\quad \alpha=0.1,
	\]
	where the  source term $f$ is given by
	\[
	\begin{array}{ll}
	f(\x,\s,t)=&(\sigma_{\tau}-\alpha-\sigma_{s})u+\pi e^{-\alpha t} s_{1}(\cos(\pi x_1)\sin(\pi x_2)\sin(\pi x_3))\\[4pt]
	&+\pi e^{-\alpha t} s_{2}(\sin(\pi x_1)\cos(\pi x_2)\sin(\pi x_3))\\[4pt]
	&+\pi e^{-\alpha t} s_{3}(\sin(\pi x_1)\sin(\pi x_2)\cos(\pi x_3)).
	\end{array}
	\]
	And the  scattering phase function   $\Phi(\s, \s')=\dfrac{ (2+2\s \cdot \s')}{4 \pi}$.
\end{example}

To verify the accuracy of the numerical approximation, 
we  discuss the discretization errors in the solution of the same example. In particular,  we have presented the discretization error in the above-defined norm  to authenticate the theoretical results. From  Table \ref{tab2}, we can see the order of converges, as expected, with the exact solution's sufficient regularity.

\begin{table}[!ht]
	\renewcommand{\arraystretch}{1}
	\caption{\label{tab2}  \it{Discretization error  for operator-splitting method of Example \ref{ex1}. }}
	{\centering
		\begin{tabular}{c c c cc }	
			\hline \hspace{2cm} &   \hspace{2cm}  & \hspace{2cm} &  \hspace{2cm}  & \hspace{2cm} \\[-6pt]
			Level  & $L^{2}$    & order & $L^{2}(0,1;L^{2}(\mathcal{G}))$    & order  \\
			\hline &&&&\\
			1  & 2.6088e-01 &  & 2.6855e-01&  \\[6pt]
			2	 &6.3417e-02  & 2.0404  & 8.1076e-02& 1.7278 \\[6pt]
			3	&  1.9910e-02 &  1.7313  &2.7095e-02& 1.5812 \\[6pt]
			4	&  8.1660e-03  & 1.2259  &9.6150e-03&1.4947\\  
			\hline
		\end{tabular}
		\par}
\end{table}

\begin{example}\label{ex3}
	Consider the model problem~\eqref{RT1} with the exact solution  
	\[
	u(\x,\s,t)=e^{-\alpha t}s_{3} \sin(\pi x_1)\sin(\pi x_2)\sin(\pi x_3),\quad \alpha=0.1, 
	\]
	where the source term $f$ is given by
	\[
	\begin{array}{ll}
	f(\x,\s,t)=&(\sigma_{\tau}-\alpha-\sigma_{s}\eta \cos \theta )u+\pi e^{-\alpha t} s_{1}s_{3}(\cos(\pi x_1)\sin(\pi x_2)\sin(\pi x_3))\\[4pt]
	&	+\pi e^{-\alpha t} s_{2}s_{3}(\sin(\pi x_1)\cos(\pi x_2)\sin(\pi x_3))\\[4pt]
	&+\pi e^{-\alpha t}s_{3}^{2}(\sin(\pi x_1)\sin(\pi x_2)\cos(\pi x_3)).
	\end{array}
	\]
	And, the  Henyey-Greenstein phase function is considered as $\Phi(\s, \s')=\dfrac{ 1}{4 \pi}\dfrac{1-\eta^{2}}{(1+\eta^{2}-2\eta \s \cdot \s')^{3/2}}$, where the anisotropy factor $\eta \in (-1,1)$. 
\end{example}	

We have discussed the convergence estimate of this test example with the anisotropy factor $\eta =0.5$. The discretization error with the convergence order is presented  in Table \ref{tab4}. These findings again confirm the error estimates of the numerical approximation achieved in the theoretical findings.

\begin{table}[!ht]
	\renewcommand{\arraystretch}{1}
	\caption{\label{tab4}  \it{Discretization error  for operator-splitting method of Example \ref{ex3}. }}
	{\centering
		\begin{tabular}{c cc c c }	
			\hline \hspace{2cm} & &  &  \hspace{2cm}  & \hspace{2cm} \\[-6pt]
			Level  & $L^{2}$    & order & $L^{2}(0,1;L^{2}(\mathcal{G}))$    & order  \\
			\hline &&&&\\
			1 & 1.8559e-01& & 1.8640e-01& \\[6pt]
			2& 7.9253e-02   &   1.2276 &8.2175e-02& 1.1816\\[6pt]
			3	& 3.9041e-02  &   1.0215 &3.9019e-02& 1.0745 \\[6pt]
			4  & 1.9017e-02&  1.0377 & 1.9381e-02 &1.0095\\  
			\hline
		\end{tabular}
		\par}
\end{table}

One can see that the convergence order is less than $1.5$, and it is due to the dependence of the angular variable on the exact solution in both the test problem. Since DG($0$) element is used, the optimal discretization error is almost first-order. The numerical results conclude that we can use tailored numerical methods in the operator-splitting finite element methods. It also explains that the convergence error is not affected by the consistency error induced by the  Lie–Trotter splitting technique in the backward Euler heterogeneous finite element method.

\section{Conclusion and Discussion}\label{sec5}
An operator-splitting finite element method for the time-dependent, high-dimensional radiative transfer equation is proposed in this paper.
The numerical scheme combines the backward Euler scheme, SUPG method, and DG($0$) for time, space, and angular discretization. The stability and consistency are established for the fully discrete scheme. Further, the convergence estimate with optimal order is derived. Moreover, the operator-splitting algorithm to compute the solution is also presented. An array of numerical experiments are performed to support the theoretical estimates and validate the proposed algorithm. The computed numerical results validate the implementation and confirm the derived error estimate.

\begin{acknowledgements}
The first author acknowledges partial support of the Science and Engineering Research Board (SERB) with the grant EMR/2016/003412.
The second author is grateful to the National Board of Higher Mathematics, Department of Atomic Energy (DAE), India, for granting postdoctoral fellowships at the Indian Institute of Science, Bangalore. 
\end{acknowledgements}

%
\section*{Declarations}
\textbf{Conflict of interest} The authors declare no competing interests.



\end{document}